\documentclass[reqno,12pt,pdftex]{amsart}

\usepackage{amsmath,amsthm,amssymb,amsfonts,amscd,mathtools}
\usepackage[all]{xy} 
\numberwithin{equation}{section}
\theoremstyle{plain} 
	\newtheorem{thm}{Theorem}[section]
	\newtheorem*{thm*}{Theorem}
	\newtheorem{cor}[thm]{Corollary}
	\newtheorem{lem}[thm]{Lemma}
	
	\newtheorem{prop}[thm]{Proposition}
	\newtheorem{conj}[thm]{Conjecture}
	\newtheorem*{conj*}{Conjecture}
	
\theoremstyle{definition}
	\newtheorem{defn}[thm]{Definition}%[section]

\theoremstyle{remark}
	\newtheorem{rem}[thm]{Remark}
	
	\newtheorem*{pf}{Proof}
\def\CC{{\mathbb C}}

\def\HH{{\mathbb H}}

\def\PP{{\mathbb P}}

\def\RR{{\mathbb R}}

\def\ZZ{{\mathbb Z}}
\def\A{{\mathcal A}}

\def\D{{\mathcal D}}
\def\E{{\mathcal E}}
\def\F{{\mathcal F}}

\def\L{{\mathcal L}}

\def\P{{\mathcal P}}

\def\R{{\mathcal R}}
\def\S{{\mathcal S}}

\def\Z{{\mathcal Z}}

\def\Aut{{\rm Aut}}

\def\Ext{{\rm Exc}}
\def\Ext{{\rm Ext}}

\def\Hom{{\rm Hom}}

\def\Sim{{\rm Sim}}
\def\Stab{{\rm Stab}}
\def\mod{{\rm mod}}

\def\ex#1{\langle #1 \rangle_{\rm ex}}
\begin{document}
\title{Full exceptional collections and stability conditions for Dynkin quivers}
\date{\today}
\author{Takumi Otani}
\address{Department of Mathematics, Graduate School of Science, Osaka University, 
Toyonaka Osaka, 560-0043, Japan}
\email{u930458f@ecs.osaka-u.ac.jp}

\maketitle
%%%%%%%%%%%%%%%%%%%%%%%%%%%%%%%%%%%%%%%%%%%%%%%%%%%%%%%%%%%%%%%%
\begin{abstract}
For a stability condition $\sigma$ on a triangulated category, Dimitrov--Katzarkov introduced the notion of a $\sigma$-exceptional collection.
In this paper, we study full $\sigma$-exceptional collections in the derived category of an acyclic quiver.
In particular, we prove that any stability condition $\sigma$ on the derived category of a Dynkin quiver admits a full $\sigma$-exceptional collection.
\end{abstract}
%%%%%%%%%%%%%%%%%%%%%%%%%%%%%%%%%%%%%%%%%%%%%%%%%%%%%%%%%%%%%%%%
\section{Introduction}
The notion of a {\em stability condition} on a triangulated category was introduced by Bridgeland.
He showed in \cite{B} that the space of stability conditions $\Stab(\D)$ on a triangulated category $\D$ has a structure of a complex manifold.
It is expected that the space of stability conditions $\Stab(\D)$ is related to various ``deformation theories''.
In order to study the space of stability conditions $\Stab(\D)$, Macr\`{i} constructed a stability condition associated with a full exceptional collection.
In several examples, such stability conditions are useful to study the homotopy type of the space of stability conditions (cf. \cite{DK3, DK4, L, M, M2}).
Motivated by Macr\`{i}'s work, Dimitrov--Katzarkov introduced the notion of a {\em $\sigma$-exceptional collection} for a stability condition $\sigma$ on a triangulated category \cite{DK}.
Roughly speaking, a $\sigma$-exceptional collection is an Ext-exceptional collection consisting of $\sigma$-stable objects such that their phases are in an interval of length one (see Definition \ref{defn : sigma-exceptional collection}).
In this paper, we study full Ext-exceptional collections in the bounded derived category $\D^b(Q)$ of finitely generated right $\CC Q$-modules of an acyclic quiver $Q$ satisfying certain conditions.
More precisely, we consider acyclic quivers with the two conditions:
\begin{enumerate}
\item[(A1)] For each $i, j = 1, \dots, \mu$, the number of arrows from $i$ to $j$ is less than one.
\item[(A2)] Let $i, k, l = 1, \dots, \mu$ such that $k < i < l$.
If there are arrows from $k$ to $i$ and from $i$ to $l$, then there are no arrows from $k$ to $l$.
\end{enumerate}
Dynkin quivers and extended Dynkin quivers except $A^{(1)}_{1,1}$ and $A^{(1)}_{1,2}$ are contained in the class of acyclic quivers satisfying (A1) and (A2).
In cases $A^{(1)}_{1,1}$ and $A^{(1)}_{1,2}$, full Ext-exceptional collections are already studied by Macr\`{i} \cite{M} and Dimitrov--Katzarkov \cite{DK}, respectively.
The following theorem is the main result of this paper.
\begin{thm}[Theorem \ref{thm : main thm}]\label{thm : main in introduction}
Let $Q$ be an acyclic quiver satisfying the conditions (A1) and (A2), and $\A$ a heart of a bounded $t$-structure on $\D^b(Q)$.
Assume that the heart $\A$ is obtained from the standard heart by iterated simple tilts.
Then, there exists a full Ext-exceptional collection $\E = (E_1, \dots, E_\mu)$ such that $\A = \ex{\E}$ and $\Sim \, \A = \{ E_1, \dots, E_\mu \}$.
\end{thm}
To prove the theorem, we use the result by King--Qiu \cite{KQ}.
They studied the relationship of simple tiltings of hearts in $\D^b(Q)$ of an acyclic quiver $Q$.
Based on simple tiltings, we construct a full $\sigma$-exceptional collection consisting of simple objects in the heart.
As a corollary of Theorem \ref{thm : main in introduction}, for a stability condition $\sigma$ whose heart is obtained from the standard heart by iterated simple tilts, there exists a full $\sigma$-exceptional collection (Corollary \ref{cor : full sigma-exceptional collection}).
For a Dynkin quiver $\vec{\Delta}$, it was shown by \cite{KV} (cf. \cite{Q}) that any heart of a bounded $t$-structure on $\D^b(\vec{\Delta})$ can be obtained from the standard heart by iterated simple tilts.
As a direct consequence, we obtain the following theorem.
\begin{thm}[Corollary \ref{cor : Dynkin}]\label{thm : Dynkin}
Let $\vec{\Delta}$ be one of the Dynkin quivers.
For each stability condition $\sigma = (Z, \P)$ on $\D^b(\vec{\Delta})$, there exists a full $\sigma$-exceptional collection.
\end{thm}
This theorem gives an affirmative answer to a conjecture by Dimitrov--Katzarkov \cite[Conjecture 7.1]{DK2}.
As a corollary of Theorem \ref{thm : Dynkin}, the space of stability conditions $\Stab(\D^b(\vec{\Delta}))$ can be described by full Ext-exceptional collections (Corollary \ref{cor : main}).
For the affine $A^{(1)}_{1,1}$ quiver and the affine $A^{(1)}_{1,2}$ quiver, it is known that for each stability condition $\sigma$ there exists a full $\sigma$-exceptional collection \cite{M, DK}.
Based on the results and Theorem \ref{thm : Dynkin}, in the case of extended Dynkin quivers the following conjecture is expected.
\begin{conj}[Conjecture \ref{conj : extended Dynkin}]
Let $Q$ be one of the extended Dynkin quivers.
For each stability condition $\sigma$ on $\D^b(Q)$, there exists $s \in \CC$ such that the heart $\P(({\rm Re}(s), {\rm Re}(s) + 1])$ can be obtained from the standard heart by iterated simple tilts.
In particular, any stability condition $\sigma$ on $\D^b(Q)$ admits a full $\sigma$-exceptional collection.
\end{conj}
\bigskip
This paper is organized as follows.
In Section \ref{sec : preliminaries}, we recall some notations of triangulated category and collect some facts related to stability conditions.
In Section \ref{sec : main results}, after recalling the notion of $\sigma$-exceptional collection for a stability condition $\sigma$ on a triangulated category, we state our main results.
Section \ref{sec : proof of thm} is devoted to proving Theorem \ref{thm : main thm}.
\bigskip
\noindent{\bf Acknowledgements.}
I am deeply grateful to my supervisor Professor Atsushi Takahashi for careful reading of this manuscript and valuable comments, and Professor Akishi Ikeda for discussions about stability conditions and constant encouragements.
I would like to thank Yu Qiu for valuable suggestions and comments.
I would also like to thank Yuichiro Goto, Kohei Kikuta and Yuuki Shiraishi for their helpful discussions. 
The author is supported by JST SPRING, Grant Number JPMJSP2138.
%%%%%%%%%%%%%%%%%%%%%%%%%%%%%%%%%%%%%%%%%%%%%%%%%%%%%%%%%%%%%%%%
\section{Preliminaries}\label{sec : preliminaries}
Throughout this paper, for a finite dimensional $\CC$-algebra $A$ (resp. a quiver $Q$), the bounded derived category of finitely generated right $A$-modules (resp. $\CC Q$-modules) is denoted by $\D^b(A) \coloneqq \D^b \mod (A)$ (resp. $\D^b(Q) \coloneqq \D^b \mod (\CC Q)$).
For convenience, we always assume that for an acyclic quiver $Q$ there are no arrows from the vertex $i$ to the other vertex $j$ when $i > j$.

Let $\D$ be a $\CC$-linear triangulated category and $[1]$ the shift functor.
For convenience, for $E, F \in \D$ we put 
\[
\Hom_\D^\bullet (E, F) \coloneqq \bigoplus_{p \in \ZZ} \Hom_\D^p (E, F)[-p] , \quad \Hom_\D^p (E, F) \coloneqq \Hom_\D (E, F[p]).
\]
A $\CC$-linear triangulated category $\D$ is said to be {\em of finite type} if the dimension of the $\CC$-vector space $\bigoplus_{p \in \ZZ} \Hom_\D^p (E, F)$ is finite for all $E, F \in \D$.

A full subcategory $\A$ of a triangulated category $\D$ will be called {\em extension closed} if whenever $E \to F \to G$ is a triangle in $\D$ with $E \in \A$ and $G \in \A$, then $F \in \A$. 
The extension closed subcategory $\ex{\S}$ of $\D$ generated by a full subcategory $\S \subset \D$ is the smallest extension closed full subcategory of $\D$ containing $\S$.
\subsection{Exceptional collections and mutations}
In this section, we recall related notions of an exceptional collection.
\begin{defn}
Let $\D$ be a $\CC$-linear triangulated category of finite type.
\begin{enumerate}
\item An object $E \in \D$ is called {\em exceptional} if $\Hom_\D(E,E) \cong \CC$ and $\Hom^p_\D(E,E) \cong 0$ when $p \ne 0$.
\item An ordered set $\E = (E_1, \dots, E_\mu)$ consisting of exceptional objects $E_1, \dots, E_\mu$ is called {\em exceptional collection} if $\Hom^p_\D(E_i, E_j) \cong 0$ for all $p \in \ZZ$ and $i > j$.
\item An exceptional collection $\E = (E_1, \dots, E_\mu)$ is called {\em Ext} if
\[
\Hom_\D^p (E_i, E_j) \cong 0,\quad  i \ne j, ~p \le 0.
\]
\item An exceptional collection $\E$ is called {\em full} if the smallest full triangulated subcategory of $\D$ containing all elements in $\E$ is equivalent to $\D$ as a triangulated category.
\item An exceptional collection $\E = (E_1, \dots, E_\mu)$ is called {\em monochromatic} if for any $i, j = 1, \dots, \mu$ such that $\Hom^\bullet_\D (E_i, E_j) \not\cong 0$, the $\ZZ$-graded $\CC$-vector space $\Hom^\bullet_\D (E_i, E_j)$ is concentrated in a single non-negative degree.
\end{enumerate}
\end{defn}
Let $\E = (E_1, \dots, E_\mu)$ be a monochromatic full Ext-exceptional collection in $\D$.
For each $i, j = 1, \dots, \mu$ such that $\Hom^\bullet_\D (E_i, E_j) \not\cong 0$, the degree in which $\Hom^\bullet_\D (E_i, E_j)$ is concentrated will be denoted by $\deg \, \Hom^\bullet_\D (E_i, E_j) \in \ZZ_{\ge 0}$.
\begin{rem}
We do not consider the degree if $\Hom^\bullet_\D (E_i, E_j)$ vanishes.
\end{rem}
For an exceptional object $E \in \D$, the smallest triangulated category containing $E$ can be identified with the bounded derived category $\D^b(\CC)$ of the $\CC$-algebra $\CC$.
Hence, the tensor functor $- \otimes E \colon \D^b(\CC) \longrightarrow \D$ is naturally defined as follows:
\[
V^\bullet \otimes E \coloneqq \bigoplus_{p \in \ZZ} V^p \otimes E [-p], \quad V^p \otimes E \coloneqq E^{\oplus \dim_\CC V^p},
\]
where $V^\bullet = \bigoplus_{p \in \ZZ} V^p [-p] \in \D^b (\CC)$.
\begin{defn}
Let $(E,F)$ be an exceptional collection.
Define two objects $\R_F E$ and $\L_E F$ by the following exact triangles respectively:
\begin{equation*}
\R_F E \longrightarrow E \overset{{\rm ev}^*}{\longrightarrow} \Hom^\bullet_\D(E, F)^* \otimes F,
\end{equation*}
\begin{equation*}
\Hom^\bullet_\D(E, F) \otimes E \overset{\rm ev}{\longrightarrow} F \longrightarrow \L_E F,
\end{equation*}
where $(-)^*$ denotes the duality $\Hom_\CC(-, \CC)$.
The object $\R_F E$ (resp. $\L_E F$) is called the {\em right mutation} of $E$ through $F$ (resp. {\em left mutation} of $F$ through $E$).
Then, $(F, \R_F E)$ and $(\L_E F, E)$ form new exceptional collections.
\end{defn}
\begin{rem}
Our definition of mutations differs from the usual one (cf. \cite[Section 1]{BP}).
In our notation, the usual left and right mutations are given by $\R_F E [1]$ and $\L_E F [-1]$, respectively.
\end{rem}
The Artin's {\it braid group} ${\rm Br}_{\mu}$ on $\mu$-stands is a group presented by the following generators and relations: 
\begin{description}
\item[{\bf Generators}] $\{b_i~|~i=1,\dots, \mu-1\}$
\item[{\bf Relations}] $b_{i}b_{j}=b_{j}b_{i}$ for $|i-j|\ge 2$, $b_{i}b_{i+1}b_{i}=b_{i+1}b_{i}b_{i+1}$ for $i=1,\dots, \mu-2$.
\end{description}
Consider the group ${\rm Br}_{\mu} \ltimes \ZZ^{\mu}$, the semi-direct product of the braid group ${\rm Br}_{\mu}$ and the 
abelian group $\ZZ^{\mu}$, defined by the group homomorphism ${\rm Br}_{\mu} \longrightarrow {\mathfrak S}_{\mu} \longrightarrow {\rm Aut}_{\ZZ} \ZZ^{\mu}$, where the first homomorphism is $b_{i}\mapsto (i, i+1)$ and the second one is induced by the natural actions of the symmetric group ${\mathfrak S}_{\mu}$ on $\ZZ^{\mu}$. 
\begin{prop}[{cf. \cite[Proposition 2.1]{BP}}]\label{prop : braid}
The group ${\rm Br}_{\mu} \ltimes \ZZ^{\mu}$ acts on the set of isomorphism classes of full exceptional collections in $\D$ by mutations and transformations:
\begin{gather*}
b_{i} \cdot (E_{1}, \dots, E_{\mu}) \coloneqq (E_{1}, \dots, E_{i-1}, E_{i+1}, \R_{E_{i+1}} E_i, E_{i+2}, \dots, E_{\mu}), \\
b^{-1}_{i} \cdot (E_{1}, \dots, E_{\mu}) \coloneqq (E_{1}, \dots, E_{i-1}, \L_{E_{i}} E_{i+1}, E_{i}, E_{i+2},\dots, E_{\mu}), \\
e_{i} \cdot (E_{1}, \dots, E_{\mu}) \coloneqq (E_{1}, \dots, E_{i-1}, E_{i} [1], E_{i+1}, \dots, E_{\mu}),
\end{gather*}
where $e_{i}$ is the $i$-th generator of $\ZZ^{\mu}$.
\end{prop}
\begin{rem}\label{rem : Ext}
For any full exceptional collection $(E_1, \dots, E_\mu)$, since $\D$ is of finite type one can choose integers $p_1, \dots, p_\mu \in \ZZ$ so that the shifted full exceptional collection $(E_1 [p_1], \dots, E_\mu [p_\mu])$ is Ext.
\end{rem}
\subsection{Hearts of bounded $t$-structures}
In this section, we collect some facts about a heart of a bounded $t$-structure.
Let $\D$ be a $\CC$-linear triangulated category of finite type.
Recall that a $t$-structure on $\D$ is a full subcategory $\F \subset \D$, satisfying $\F[1] \subset \F$, such that for every object $E \in \D$, there is a triangle $F \to E \to G$ in $\D$ with $F \in \F$ and $G \in \F^\perp$, where $\F^\perp$ is a full subcategory given by 
\[
\F^\perp \coloneqq \{ G \in \D \mid \Hom_\D (F, G) = 0, ~ F \in \F \}.
\]
A $t$-structure $\P$ is said to be {\em bounded} if 
\[
\D = \bigcup_{i,j \in \ZZ} \F [i] \cap \F^\perp [j].
\]
A {\em heart} $\A$ of a $t$-structure $\F$ is the full subcategory $\A \coloneqq \F \cap \F^\perp [1]$.
It was proved in \cite{BBD} that a heart $\A$ of a $t$-structure is an abelian category, with the short exact sequences in $\A$ being precisely the triangles in $\D$ all of whose vertices are objects of $\A$.
For simplicity, we use the term {\em heart} to mean the heart of a bounded $t$-structure.
Let $\A$ be a heart in $\D$ and $S \in \A$ a simple object.
Define full subcategories
\[
{}^\perp S \coloneqq \{ E \in \A \mid \Hom_\A(E, S) = 0 \}, \quad
S^\perp \coloneqq \{ E \in \A \mid \Hom_\A(S, E) = 0 \}.
\]
Then, one can consider two full subcategories 
\[
\A_S^\sharp \coloneqq \ex{S[1], {}^\perp S}, \quad 
\A_S^\flat \coloneqq \ex{S^\perp, S[-1]}.
\]
It is known that the subcategories $\A_S^\sharp$ and $\A_S^\flat$ are hearts in $\D$.
The heart $\A_S^\sharp$ (resp. $\A_S^\flat$) is called the {\em forward simple tilt} (resp. {\em backward simple tilt}) of $\A$ by $S$.
Denote by $\Sim \, \A$ a complete set of simple modules in a heart $\A$.
\begin{prop}[{\cite[Proposition 5.4]{KQ}}]\label{prop : Exchange graph}
Let $S$ be a simple object in $\A$ satisfying $\Ext^1(S, S) \cong 0$.
Assume that $\#\Sim \, \A$ is finite and $\A = \ex{\Sim \, \A}$.
Then, we have 
\begin{eqnarray*}
\Sim \, \A_S^\sharp & = & \{ S[1] \} \cup \{ \psi^\sharp_{S}(T) \mid T \in \Sim \, \A, ~ T \not\cong S \}, \\
\Sim \, \A_S^\flat & = & \{ S[-1] \} \cup \{ \psi^\flat_{S}(T) \mid T \in \Sim \, \A, ~ T \not\cong S \},
\end{eqnarray*}
where $\psi_{S}^\sharp(T)$ and $\psi_S^\flat(T)$ are defined by the following exact triangles:
\begin{subequations}\label{eq : twist}
\begin{equation}\label{eq : twist sharp}
\psi_{S}^\sharp(T) \longrightarrow T \longrightarrow \Hom^1(T, S)^* \otimes S[1],
\end{equation}
\begin{equation}\label{eq : twist flat}
\Hom^1(S, T) \otimes S[-1] \longrightarrow T \longrightarrow \psi_{S}^\flat(T).
\end{equation}
\end{subequations}
and $\A_S^\sharp = \ex{\Sim \, \A_S^\sharp}$ and $\A_S^\flat = \ex{\Sim \, \A_S^\flat}$.
\qed
\end{prop}
For an acyclic quiver, King--Qiu also showed the following
\begin{prop}[{\cite[Proposition 6.4]{KQ}}]\label{prop : strongly monochromaticity}
Let $Q$ be an acyclic quiver and $\A$ a heart in $\D^b(Q)$.
For any distinct simple objects $S$ and $T$ in the heart $\A$, $\Hom^\bullet_{\D^b(Q)} (S, T)$ and $\Hom^\bullet_{\D^b(Q)} (T, S)$ are concentrated in a single positive degree.
Moreover, we have $\Hom^\bullet_{\D^b(Q)} (S, T) \cong 0$ or $\Hom^\bullet_{\D^b(Q)} (T, S) \cong 0$.
\qed
\end{prop}
In the Dynkin case, the following result is known.
It plays an important role in this paper.
\begin{prop}[{\cite{KV} (cf. \cite[Appendix A]{Q})}]\label{prop : Dynkin Exchange graph}
Any heart in $\D^b(\vec{\Delta})$ can be obtained from the standard heart $\mod(\CC \vec{\Delta})$ by iterated simple tilts.
\qed
\end{prop}
\subsection{Stability conditions}
In this section, we recall the notion of a stability condition on a triangulated category.
Let $\D$ be a $\CC$-linear triangulated category of finite type.
Denote by $K_0(\D)$ the Grothendieck group of the triangulated category $\D$.
We will denote by $\mu \in \ZZ_{\ge 0}$ the rank of the Grothendieck group $K_0(\D)$.
\begin{defn}[{\cite[Definition 1.1]{B}}]\label{defn : stability condition}
A {\em stability condition} on $\D$ consists of a group homomorphism $Z \colon K_0(\D)\longrightarrow\CC$, which is called {\em central charge}, and the family of full additive subcategories $\P = \{ \P(\phi) \}_{\phi \in \RR}$, called {\em slicing}, satisfying the following axioms:
\begin{enumerate}
\item if $E \in \P(\phi)$ then $Z(E) = m(E) \exp(\sqrt{-1} \pi \phi)$ for some $m(E) \in \RR_{>0}$.
\item for all $\phi \in \RR$, $\P (\phi + 1) = \P(\phi) [1]$,
\item if $\phi_1 > \phi_2$ and $A_i \in \P(\phi_i)$ then $\Hom_\D (A_1, A_2) = 0$, 
\item for each nonzero object $E \in \D$ there exists a finite sequence of real numbers
\[
\phi_1 > \phi_2 > \dots > \phi_n
\]
and a collection of triangles
\small
\[
\xymatrix{
0 = F_0 \ar[rr] & & F_1 \ar[rr] \ar[ld] & & F_2 \ar[r] \ar[ld] & \cdots \ar[r] & F_{n-1} \ar[rr] & & F_n = E \ar[ld] \\
& A_1 \ar@{-->}[ul] & & A_2 \ar@{-->}[ul] & & & & A_n \ar@{-->}[ul] & 
}
\]
\normalsize
with $A_i \in \P(\phi_i)$ for all $i = 1, \dots, n$.
\item (support property) There exists a constant $C > 0$ such that for all nonzero object $E \in \P(\phi)$ for some $\phi \in \RR$, we have 
\[
\| E \| < C | Z(E) |
\]
where $\| \cdot \|$ denotes a norm on $K_0(\D) \otimes_\ZZ \RR$.
\end{enumerate}
\end{defn}
A nonzero object $E \in \P(\phi)$ is called {\em $\sigma$-semistable} of phase $\phi$, and a simple object of $E \in \P(\phi)$ is called {\em $\sigma$-stable}.
For any interval $I \subset \RR$, we put $\P(I) \coloneqq \ex{\P(\phi) \mid \phi \in I}$.
Then, the full subcategory $\P((0, 1])$ is a heart in $\D$, hence an abelian category.
Conversely, a stability condition can be described as a heart and a certain function on the heart.
First, we recall the notion of stability function.
We denote by $K_0(\A)$ the Grothendieck group of an abelian category $\A$.
\begin{defn}[{\cite[Definition 2.1]{B}}]
Let $\A$ be a heart in $\D$.
A {\em stability function} on $\A$ is a group homomorphism $Z \colon K_0(\A) \longrightarrow \CC$ such that for all nonzero object $E \in \A$ the complex number $Z(E)$ lies in the semiclosed upper half plane 
$\HH_- \coloneqq \{ r e^{\sqrt{-1} \pi \phi} \in \CC \mid r > 0, ~ 0 < \phi \le 1 \}$.
\end{defn}
Given a stability function $Z \colon K_0(\A) \longrightarrow \CC$, the {\em phase} of a nonzero object $E \in \A$ is defined to be the real number $\phi(E) \coloneqq (1 / \pi) {\rm arg} \, Z(E) \in (0, 1]$.
A nonzero object $E \in \A$ is {\em semistable} (resp. {\em stable}) if we have $\phi(A) \le \phi(E)$ (resp. $\phi(A) < \phi(E)$) for all nonzero subobjects $A \subset E$.
We say that a stability function $Z \colon K_0(\A) \longrightarrow \CC$ satisfies the {\em Harder--Narasimhan property} if each nonzero object $E \in \A$ admits a filtration 
\[
0 = F_0 \subset F_1 \subset \cdots \subset F_{n - 1} \subset F_n = E
\]
such that $F_i / F_{i - 1}$ is semistable for $i = 1, \dots, n$ with $\phi(F_1/F_0) > \phi(F_2/F_1) > \dots > \phi(F_n/F_{n - 1})$.
We also say that a stability function  $Z \colon K_0(\A) \longrightarrow \CC$ satisfies the {\em support property} if there exists a constant $C > 0$ such that for all semistable objects $E \in \A$, we have $\| E \| < C | Z(E) |$.
\begin{prop}[{\cite[Proposition 5.3]{B}}]
To give a stability condition on a triangulated category $\D$ is equivalent to giving a bounded $t$-structure on $\D$ and a stability function on its heart with the Harder-Narasimhan property and the support property.
\qed
\end{prop}
Denote by $\Stab(\D)$ the space of stability conditions on $\D$.
\begin{prop}[{\cite[Theorem 1.2]{B}}]\label{prop : bridgeland}
There exists a natural topology on $\Stab(\D)$ such that the forgetful morphism 
\begin{equation*}
\Z \colon \Stab(\D) \longrightarrow \Hom_\ZZ(K_0(\D), \CC), \quad 
(Z, \P) \mapsto Z
\end{equation*}
is a local homeomorphism.
In particular, the space $\Stab(\D)$ has a structure of complex manifolds of dimension $\mu$.
\qed
\end{prop}
There are natural group actions on $\Stab(\D)$ commuting with each other.
The first one is a $\CC$-action 
\begin{equation*}
s \cdot (Z, \P) = (e^{- \pi \sqrt{-1} s} \cdot Z, \P_{{\rm Re} (s)}), \quad s \in \CC, 
\end{equation*}
where $\P_{{\rm Re} (s)}(\phi) \coloneqq \P (\phi + {\rm Re} (s))$.
The other action is given by the group of autoequivalences $\Aut(\D)$
\begin{equation*}
\Phi \cdot (Z, \P) = (Z \circ F^{-1}, \Phi(\P)), \quad \Phi \in \Aut(\D).
\end{equation*}
For a heart $\A$ in $\D$, we can consider the subset $U(\A) \subset \Stab(\D)$ defined by 
\[
U(\A) \coloneqq \{ (Z, \P) \in \Stab(\D) \mid \P((0, 1]) = \A \}.
\]
A heart $\A$ is said to be {\em of finite length} if $\A$ is artinian and noetherian.
When the heart $\A$ is of finite length and has finitely many simple objects, any stability function $Z \colon K_0(\A) \longrightarrow \CC$ satisfies the Harder--Narasimhan property.
Moreover, due to \cite[Proposition B.4]{BM}, a stability function $Z \colon K_0(\A) \longrightarrow \CC$ on a heart $\A$ of finite length with finitely many simple objects satisfies the support property.
Therefore, we have the following
\begin{prop}[{\cite[Lemma 5.2]{B2}}]\label{prop : finite length heart}
Let $\A$ be a heart of finite length in $\D$ with $\Sim \, \A = \{ S_1, \dots, S_\mu\}$.
Then, we have an isomorphism 
\[
U(\A) \overset{\cong}{\longrightarrow} \HH_-^\mu, \quad 
(Z, \P) \mapsto (Z(S_1), \dots, Z(S_\mu)).
\]
\qed
\end{prop}
%%%%%%%%%%%%%%%%%%%%%%%%%%%%%%%%%%%%%%%%%%%%%%%%%%%%%%%%%%%%%%%%
\section{Main results}\label{sec : main results}
In this section, we study a relation between exceptional collections and stability conditions.
The following proposition proved by Macr\`{i} associates full exceptional collections to stability conditions.
\begin{prop}[{\cite[Lemma 3.14 and Lemma 3.16]{M}}]\label{prop : Ext and heart}
Let $\E = (E_1, \dots, E_\mu)$ be a full Ext-exceptional collection in $\D$.
Then, the extension closed subcategory $\ex{\E}$ is a heart of finite length in $\D$ and $\Sim \, \ex{\E} = \{ E_1, \dots, E_\mu \}$.
\qed
\end{prop}
If $\E$ is a full Ext-exceptional collection, then we can construct a stability condition $\sigma = (Z, \P)$ such that $\P((0, 1]) = \ex{\E}$ by Proposition \ref{prop : finite length heart}.
Due to Remark \ref{rem : Ext}, we have the following
\begin{cor}[{\cite[Section 3.3]{M}}]\label{cor : full exceptional collection}
For a full exceptional collection $\E = (E_1, \dots, E_\mu)$ in $\D$, there exists a stability condition $\sigma$ on $\D$ such that $E_i$ is $\sigma$-stable for each $i = 1, \dots, \mu$.
\qed
\end{cor}
Motivated by Macr\`{i}'s work, Dimitrov--Katzarkov introduced the notion of a $\sigma$-exceptional collection for a stability condition $\sigma$.
\begin{defn}[{\cite[Definition 3.17]{DK}}]\label{defn : sigma-exceptional collection}
Let $\sigma = (Z, \P) \in \Stab(\D)$ be a stability condition on $\D$.
An exceptional collection $\E = (E_1, \dots, E_\mu)$ in $\D$ is called {\em $\sigma$-exceptional collection} if the following three properties hold:
\begin{enumerate}
\item For each $i = 1, \dots, \mu$, the object $E_i$ is $\sigma$-semistable.
\item $\E$ is an Ext-exceptional collection.
\item There exists a real number $r \in \RR$ such that $r < \phi(E_i) \le r + 1$ for $i = 1, \dots, \mu$.
\end{enumerate}
\end{defn}
By definition, a full Ext-exceptional collection $\E$ is a full $\sigma$-exceptional collection of a stability condition $\sigma$ given in Proposition \ref{prop : Ext and heart}.

We shall prepare to state the main results.
For an acyclic quiver $Q$, we introduce the following two conditions:
\begin{enumerate}
\item[(A1)] For each $i, j = 1, \dots, \mu$, the number of arrows from $i$ to $j$ is less than one.
\item[(A2)] Let $i, k, l = 1, \dots, \mu$ such that $k < i < l$.
If there are arrows from $k$ to $i$ and from $i$ to $l$, then there are no arrows from $k$ to $l$.
\end{enumerate}
Dynkin quivers and extended Dynkin quivers except $A^{(1)}_{1,1}$ and $A^{(1)}_{1,2}$ satisfy the above two conditions (A1) and (A2).
The following theorem is the main result of this paper.
\begin{thm}\label{thm : main thm}
Let $Q$ be an acyclic quiver satisfying the conditions (A1) and (A2), and $\A$ a heart in $\D^b(Q)$.
Assume that the heart $\A$ is obtained from the standard heart by iterated simple tilts.
Then, there exists a monochromatic full Ext-exceptional collection $\E = (E_1, \dots, E_\mu)$ such that $\A = \ex{\E}$ and $\Sim \, \A = \{ E_1, \dots, E_\mu \}$.
\end{thm}
We prove Theorem \ref{thm : main thm} in Section \ref{sec : proof of thm}.
As a consequence of this theorem, we obtain the following result.
\begin{cor}\label{cor : full sigma-exceptional collection}
Let $Q$ be an acyclic quiver satisfying the conditions (A1) and (A2), and $\sigma = (Z, \P)$ a stability condition on $\D^b(Q)$. 
Assume that there exists $s \in \CC$ such that the heart $\P_{{\rm Re}(s)}((0, 1])$ of the stability condition $s \cdot \sigma$ is obtained from the standard heart by iterated simple tilts.
Then, there exists a monochromatic full $\sigma$-exceptional collection.
\end{cor}
\begin{pf}
The statement follows from Theorem \ref{thm : main thm}.
\qed
\end{pf}
In particular, when $Q$ is a Dynkin quiver $\vec{\Delta}$, we have the following 
\begin{cor}\label{cor : Dynkin}
Let $\vec{\Delta}$ be one of the Dynkin quivers.
For each stability condition $\sigma = (Z, \P) \in \Stab(\D^b(\vec{\Delta}))$, there exists a monochromatic full $\sigma$-exceptional collection $\E = (E_1, \dots, E_\mu)$. 
Moreover, we can choose $\E = (E_1, \dots, E_\mu)$ so that $\P((0, 1]) = \ex{\E}$ and $\Sim \, \P((0, 1]) = \{ E_1, \dots, E_\mu \}$.
\end{cor}
\begin{pf}
Due to Proposition \ref{prop : Dynkin Exchange graph}, any heart in $\D^b(\vec{\Delta})$ can be obtained from the standard heart by iterated simple tilts.
Therefore, one can obtain the statement by Theorem \ref{cor : full sigma-exceptional collection}.
\qed
\end{pf}
The corollary solves the conjecture given by Dimitrov--Katzarkov (\cite[Conjecture 7.1]{DK2}).
As a corollary, any stability condition can be described by a full $\sigma$-exceptional collection and a stability function.
Namely, we have the following
\begin{cor}[{cf. \cite[Theorem 2.12]{W}}]\label{cor : main}
We have 
\[
\Stab(\D^b(\vec{\Delta})) \cong \bigcup \{ U(\ex{\E}) \mid \text{$\E$ is a monochromatic full Ext-exceptional collection} \}.
\]
\end{cor}
\begin{pf}
The statement follows from Proposition \ref{prop : Ext and heart} and Corollary \ref{cor : Dynkin}.
\qed
\end{pf}
It is expected that the space of stability conditions $\Stab(\D^b(\vec{\Delta}))$ can be identified with the universal unfolding (deformation) space of a polynomial $f \colon \CC^3 \longrightarrow \CC$ of type ADE (cf. \cite{B2, BQS, HKK, T}).
We hope that Corollary \ref{cor : main} is useful to solve the problem.
In several cases of extended Dynkin quivers, the existence of full $\sigma$-exceptional collections is known.
\begin{prop}[{\cite[Section 4]{M} for $Q = A^{(1)}_{1,1}$, \cite[Theorem 1.1]{DK} for $Q = A^{(1)}_{1,2}$}]\label{prop : ext Dynkin}
Let $Q$ be the affine $A^{(1)}_{1,1}$ quiver or the affine $A^{(1)}_{1,2}$ quiver.
For each stability condition $\sigma$ on $\D^b(Q)$, there exists a full $\sigma$-exceptional collection.
\qed
\end{prop}
\begin{rem}
Okada classified hearts in $\D^b(A_{1, 1}^{(1)})$ on which we can impose stability function with the Harder--Narasimhan property \cite[Corollary 3.4]{O}.
Based on this result, one can check Conjecture \ref{conj : extended Dynkin} directly.
\end{rem}
Motivated by Corollary \ref{cor : Dynkin} and Proposition \ref{prop : ext Dynkin}, the following conjecture is expected.
\begin{conj}\label{conj : extended Dynkin}
Let $Q$ be one of the extended Dynkin quivers.
For each stability condition $\sigma$ on $\D^b(Q)$, there exists $s \in \CC$ such that the heart $\P(({\rm Re}(s), {\rm Re}(s) + 1])$ can be obtained from the standard heart by iterated simple tilts.
In particular, any stability condition $\sigma$ on $\D^b(Q)$ admits a full $\sigma$-exceptional collection.
\end{conj}
\begin{rem}
When $Q$ is the affine $A^{(1)}_{p, q}$ quiver, Haiden--Katzarkov--Kontsevich proved that $\Stab(\D^b(Q))$ is identified with the moduli space of framed exponential type quadratic differentials on annulus with marked points and a grading \cite[Theorem 6.2]{HKK}.
This result implies that any stability condition has a heart of finite length with finitely many simple objects up to the $\CC$-action.
Hence, in order to prove Conjecture \ref{conj : extended Dynkin} for the affine $A^{(1)}_{p, q}$ quiver, it is sufficient to show that any heart of finite length with finitely many simple objects can be obtained from the standard heart by iterated simple tilts.
\end{rem}
%%%%%%%%%%%%%%%%%%%%%%%%%%%%%%%%%%%%%%%%%%%%%%%%%%%%%%%%%%%%%%%%
\section{Proof of Theorem \ref{thm : main thm}}\label{sec : proof of thm}
The sketch of the proof is as follows.
Firstly, we show the statement in the case of the standard heart $\A = \mod \, (\CC Q)$ (Proposition \ref{prop : std heart}).
Secondly, we show the general case.
We prove that a monochromatic full Ext-exceptional collection in a heart induces a new one in the simple tilted heart (Proposition \ref{prop : tilting and full exceptional collection}).
\bigskip
For a monochromatic full Ext-exceptional collection $\E = (E_1, \dots, E_\mu)$, it will be convenient to write $p_{i, j} = \deg \, \Hom^\bullet_\D (E_i, E_j)$ if $\Hom^\bullet_\D (E_i, E_j) \not\cong 0$.
In order to show the existence of a full $\sigma$-exceptional collection for a given stability condition $\sigma$ on $\D$, we introduce the following two conditions (E1) and (E2) for monochromatic full Ext-exceptional collections:
\begin{itemize}
\item[(E1)] For each $i, j = 1, \dots, \mu$ we have $\dim_\CC \bigoplus_{p \in \ZZ} \Hom^p_\D (E_i, E_j) \le 1$.
\item[(E2)] Let $i ,k ,l = 1,\dots, \mu$ such that $k < i < l$.
If $\Hom^\bullet_\D (E_k, E_l) \not\cong 0$, $\Hom^\bullet_\D (E_k, E_i) \not\cong 0$ and $\Hom^\bullet_\D (E_i, E_l) \not\cong 0$, then we have 
\begin{equation*}
p_{k, l} = p_{k, i} + p_{i, l}.
\end{equation*}
\end{itemize}
\begin{prop}\label{prop : std heart}
Let $Q$ be an acyclic quiver satisfying the condition (A1) and (A2).
There exists a monochromatic full Ext-exceptional collection $\E = (E_1, \dots, E_\mu)$ in $\D^b(Q)$ such that $\Sim \, \mod(\CC Q) = \{ E_1, \dots, E_\mu \}$ and $\E$ satisfies the conditions (E1) and (E2).
\end{prop}
\begin{pf}
The simple $\CC Q$-module $S_i$ corresponding to the vertex $i$ is an exceptional object in $\D$.
It is known that the exceptional collection $\E \coloneqq (S_1, \dots, S_\mu)$ consisting of simple $\CC Q$-modules in the standard heart $\mod(\CC Q)$ is full.
Moreover, since $\dim_\CC \Ext^1_{\CC Q} (S_i, S_j)$ is equal to the number of arrows from the vertex $i$ to $j$, the conditions (A1) and (A2) imply (E1) and (E2), respectively.
\qed
\end{pf}
Let $\E = (E_1, \dots, E_\mu)$ a monochromatic full Ext-exceptional collection in $\D$.
Fix $i = 1,\dots, \mu$.
We can assume, by reordering $\E$ if necessary, that if $\Hom_\D^1 (E_j, E_i) \not\cong 0$ and $\Hom_\D^\bullet (E_{j'}, E_i) \not\cong 0$ for $j < j' < i$ then we have $p_{j', i} = 1$.
Indeed, if $\Hom_\D^1 (E_j, E_i) \not\cong 0$, $\Hom_\D^\bullet (E_{j'}, E_i) \not\cong 0$ and $p_{j', i} \ge 2$, then the condition (E2) implies $\Hom^\bullet_\D (E_{j'}, E_j) \cong 0$.
Define two indices $i^\sharp$ and $i^\flat$ by
\begin{eqnarray*}
i^\sharp & \coloneqq & \min \{ j \in \{ 1, \dots, i \} \mid \Hom^\bullet_\D (E_j, E_i) \not\cong 0, ~ p_{j, i} \le 1 \}, \\
i^\flat & \coloneqq & \max \{ j \in \{ i , \dots, \mu \} \mid \Hom^\bullet_\D (E_i, E_j) \not\cong 0, ~ p_{i, j} \le 1 \}.
\end{eqnarray*}
We also define the ordered sets $\E^\sharp_i$ and $\E^\flat_i$ by 
\begin{eqnarray*}
\E^\sharp_i & \coloneqq & ( E_1, \dots, E_{i^\sharp - 1}, E_i [1], \psi^\sharp_{E_i} (E_{i^\sharp}), \dots, \psi^\sharp_{E_i} (E_{i - 1}), E_{i + 1}, \dots, E_\mu ), \\
\E^\flat_i & \coloneqq & ( E_1, \dots, E_{i - 1},\psi^\flat_{E_i} (E_{i + 1}), \dots, \psi^\flat_{E_i} (E_{i^\flat}), E_i [-1], E_{i^\flat + 1}, \dots, E_\mu ),
\end{eqnarray*}
where $\psi^\sharp_{E_i} (E_j)$ and $\psi^\flat_{E_i} (E_j)$ are objects defined in \eqref{eq : twist}.
\begin{prop}\label{prop : tilting and full exceptional collection}
Let $Q$ be an acyclic quiver and $\A$ a heart in $\D = \D^b(Q)$. 
Suppose that there exists a monochromatic full Ext-exceptional collection $\E = (E_1, \dots, E_\mu)$ such that $\A = \ex{\E}$, $\Sim \, \A = \{ E_1, \dots, E_\mu \}$ and $\E$ satisfies the conditions (E1) and (E2).
Then, we have the following:
\begin{subequations}\label{eq : Sim E}
\begin{enumerate}
\item For each $i = 1,\dots, \mu$, the ordered set $\E^\sharp_i$ is a monochromatic full Ext-exceptional collection and satisfies the conditions (E1) and (E2). 
Moreover, we have $\A^\sharp_{E_i} = \ex{\E^\sharp_i}$ and 
\begin{equation}\label{eq : Sim E sharp}
\Sim \, \A^\sharp_{E_i} = \{ E_1, \dots, E_{i^\sharp - 1}, E_i [1], \psi^\sharp_{E_i} (E_{i^\sharp}), \dots, \psi^\sharp_{E_i} (E_{i - 1}), E_{i + 1}, \dots, E_\mu \}.
\end{equation}
\item For each $i = 1,\dots, \mu$, the ordered set $\E^\flat_i$ is a monochromatic full Ext-exceptional collection and satisfies the conditions (E1) and (E2). 
Moreover, we have $\A^\flat_{E_i} = \ex{\E^\flat_i}$ and 
\begin{equation}\label{eq : Sim E flat}
\Sim \, \A^\flat_{E_i} = \{ E_1, \dots, E_{i - 1},\psi^\flat_{E_i} (E_{i + 1}), \dots, \psi^\flat_{E_i} (E_{i^\flat}), E_i [-1], E_{i^\flat + 1}, \dots, E_\mu \}.
\end{equation}
\end{enumerate}
\end{subequations}
\end{prop}
\begin{pf}
We only prove the first statement (i).
However, a very similar proof works for the second statement (ii).
By assumption, for each $j = i^\sharp, \dots, i - 1$ we have 
\[
\Hom^\bullet_\D (E_j, E_i)^* \otimes E_i \cong (\Hom^1_\D (E_j, E_i) [-1])^* \otimes E_i \cong \Hom^1_\D(E_j, E_i)^* \otimes E_i [1].
\]
This identification yields an isomorphism $\psi^\sharp_{E_i} (E_j) \cong \R_{E_i} E_j$ for $j = i^\sharp, \dots, i - 1$.
Proposition \ref{prop : braid} implies that the ordered set $\E^\sharp_i$ forms a full exceptional collection.

In order to see that the full exceptional collection $\E^\sharp_i$ is monochromatic and Ext, we split the proof into four lemmas:
\begin{lem}\label{lem : Hom 2 sharp}
Let $k = 1, \dots, i^\sharp - 1$ and $j = i^\sharp, \dots, i - 1$.
\begin{subequations}
\begin{enumerate}
\item If $\Hom^1_\D (E_j, E_i) \cong 0$ or $\Hom^\bullet_\D (E_k, E_i) \cong 0$, then we have 
\begin{equation}\label{eq : Hom 2 sharp 1}
\Hom^\bullet_\D (E_k, \psi^\sharp_{E_i} (E_j)) \cong \Hom^\bullet_\D (E_k, E_j).
\end{equation}
\item If $\Hom^1_\D (E_j, E_i) \not\cong 0$, $\Hom^\bullet_\D (E_k, E_i) \not\cong 0$ and $\Hom^\bullet_\D (E_k, E_j) \cong 0$, then we have 
\begin{equation}\label{eq : Hom 2 sharp 2}
\Hom^p_\D (E_k, \psi^\sharp_{E_i} (E_j)) \cong 
\begin{cases}
\Hom^1_\D (E_j, E_i)^* \otimes \Hom^p_\D (E_k, E_i), & p = p_{k, i}, \\
0, & p \ne p_{k, i}.
\end{cases}
\end{equation}
\item If $\Hom^1_\D (E_j, E_i) \not\cong 0$, $\Hom^\bullet_\D (E_k, E_i) \not\cong 0$ and $\Hom^\bullet_\D (E_k, E_j) \not\cong 0$, then we have 
\begin{equation}\label{eq : Hom 2 sharp 3}
\Hom^\bullet_\D (E_k, \psi^\sharp_{E_i} (E_j)) \cong 0.
\end{equation}
\end{enumerate}
\end{subequations}
\end{lem}
\begin{pf}
If $\Hom^1_\D (E_j, E_i) \cong 0$, the exact triangle $\psi^\sharp_{E_i} (E_j) \to E_j \to \Hom^1_\D(E_j, E_i)^* \otimes E_i [1]$ implies $\psi^\sharp_{E_i} (E_j) \cong E_j$.
Hence, we have \eqref{eq : Hom 2 sharp 1}.

We assume $\Hom^1_\D (E_j, E_i) \not\cong 0$.
By the exact triangle $\psi^\sharp_{E_i} (E_j) \to E_j \to \Hom^1_\D(E_j, E_i)^* \otimes E_i [1]$, we have the long exact sequence
\small
\begin{align*}
\cdots & \longrightarrow \Hom^{p - 1}_\D (E_k, E_j) \overset{\delta_{p - 1}}{\longrightarrow}  \Hom^{p - 1}_\D (E_k, \Hom^1_\D(E_j, E_i)^* \otimes E_i [1]) \longrightarrow \\
\longrightarrow \Hom^p_\D (E_k, \psi^\sharp_{E_i} (E_j)) & \longrightarrow \Hom^p_\D (E_k, E_j) \overset{\delta_p}{\longrightarrow} \Hom^p_\D (E_k, \Hom^1_\D(E_j, E_i)^* \otimes E_i [1]) \longrightarrow \\
\longrightarrow \Hom^{p + 1}_\D (E_k, \psi^\sharp_{E_i} (E_j)) & \longrightarrow \Hom^{p + 1}_\D (E_k, E_j) \overset{\delta_{p + 1}}{\longrightarrow} \cdots, 
\end{align*}
\normalsize
where $\delta_p \colon \Hom^p_\D (E_k, E_j) \longrightarrow \Hom^p_\D (E_k, \Hom^1_\D(E_j, E_i)^* \otimes E_i [1])$ denotes the $\CC$-linear map induced by the morphism $E_j \to \Hom^1_\D(E_j, E_i)^* \otimes E_i [1]$.
Note that we have
\[
\Hom^p_\D (E_k, \Hom^1_\D(E_j, E_i)^* \otimes E_i [1]) \cong \Hom^1_\D(E_j, E_i)^* \otimes \Hom^{p + 1}_\D (E_k, E_i).
\]
Hence, if we have $\Hom^\bullet_\D (E_k, E_i) \cong 0$, the long exact sequence implies the isomorphism \eqref{eq : Hom 2 sharp 1}.

We assume $\Hom^\bullet_\D (E_k, E_i) \not\cong 0$.
If $\Hom^\bullet_\D (E_k, E_j) \cong 0$, then the long exact sequence implies \eqref{eq : Hom 2 sharp 2}.

Next, we consider the last case $\Hom^\bullet_\D (E_k, E_j) \not\cong 0$.
Since the condition (E2) implies $p_{k, j} = p_{k, i} - 1$, the $\CC$-linear map $\delta_p$ is nontrivial for $p = p_{k, j}$.
By the condition (E1), the $\CC$-linear map $\delta_p$ can be identified with the map $\delta_p \colon \CC \longrightarrow \CC$.
If $\delta_p$ is the zero map, we have $\Hom^p_\D (E_k, \psi^\sharp_{E_i} (E_j)) \cong \CC$ and $\Hom^{p + 1}_\D (E_k, \psi^\sharp_{E_i} (E_j)) \cong \CC$.
However, it contradicts Proposition \ref{prop : strongly monochromaticity}.
Therefore, the $\CC$-linear map $\delta_p \colon \CC \longrightarrow \CC$ is an isomorphism and we obtain \eqref{eq : Hom 2 sharp 3}.
\qed
\end{pf}
\begin{lem}\label{lem : Hom 1 sharp}
Let $j = i^\sharp, \dots, i-1$.
We have 
\begin{equation*}
\Hom^p_\D (E_i [1], \psi^\sharp_{E_i} (E_j)) \cong 
\begin{cases}
\Hom^1_\D(E_j, E_i)^*, & p = 1 \\
0, & p \ne 1.
\end{cases}
\end{equation*}
\end{lem}
\begin{pf}
Since $\psi^\sharp_{E_i} (E_j) \cong \R_{E_i} E_j$, the statement is known.
However, we give a proof here.

By the exact triangle $\psi^\sharp_{E_i} (E_j) \to E_j \to \Hom^1_\D(E_j, E_i)^* \otimes E_i [1]$, we have the long exact sequence
\small
\begin{alignat*}{2}
& && \mathllap{\cdots} \longrightarrow \Hom^{p - 1}_\D (E_i [1], \Hom^1_\D(E_j, E_i)^* \otimes E_i [1]) \longrightarrow \\
& \longrightarrow \Hom^p_\D (E_i [1], \psi^\sharp_{E_i} (E_j)) \longrightarrow \Hom^p_\D (E_i [1], E_j) && \longrightarrow \Hom^p_\D (E_i [1], \Hom^1_\D(E_j, E_i)^* \otimes E_i [1]) \longrightarrow \\
& \longrightarrow \Hom^{p + 1}_\D (E_i [1], \psi^\sharp_{E_i} (E_j)) \longrightarrow \cdots. &&
\end{alignat*}
\normalsize
Since $\Hom^\bullet_\D (E_i, E_j) \cong 0$, for each $p \in \ZZ$ we have 
\begin{eqnarray*}
\Hom^p_\D (E_i [1], \psi^\sharp_{E_i} (E_j)) & \cong & \Hom^{p - 1}_\D (E_i [1], \Hom^1_\D(E_j, E_i)^* \otimes E_i [1]) \\
& \cong & \Hom^{p - 1}_\D (E_i, \Hom^1_\D(E_j, E_i)^* \otimes E_i) \\
& \cong & \Hom^1_\D(E_j, E_i)^* \otimes \Hom^{p - 1}_\D (E_i, E_i) .
\end{eqnarray*}
We obtain the statement since $\Hom^\bullet_\D (E_i, E_i) \cong \CC$.
\qed
\end{pf}
\begin{lem}\label{lem : Hom 4 sharp}
Let $j, j' =i^\sharp, \dots, i - 1$ such that $j < j'$.
We have 
\begin{equation*}
\Hom^\bullet_\D (\psi^\sharp_{E_i} (E_j), \psi^\sharp_{E_i} (E_{j'})) \cong \Hom^\bullet_\D (E_j, E_{j'}).
\end{equation*}
\end{lem}
\begin{pf}
Since $\psi^\sharp_{E_i} (E_j) \cong \R_{E_i} E_j$, the statement is also known.
However, we give a proof here.

By the exact triangle $\psi^\sharp_{E_i} (E_{j'}) \to E_{j'} \to \Hom^1_\D(E_{j'}, E_i)^* \otimes E_i [1]$, we have the long exact sequence
\small
\begin{alignat*}{2}
& && \mathllap{\cdots \longrightarrow \Hom^{p - 1}_\D (\psi^\sharp_{E_i} (E_j), \Hom^1_\D(E_{j'}, E_i)^* \otimes E_i [1])} \longrightarrow \\
& \longrightarrow \Hom^p_\D (\psi^\sharp_{E_i} (E_j), \psi^\sharp_{E_i} (E_{j'})) \longrightarrow \Hom^p_\D (\psi^\sharp_{E_i} (E_j), E_{j'}) \longrightarrow \Hom^p_\D (\psi^\sharp_{E_i} (E_j), \Hom^1_\D(E_{j'}, E_i)^* \otimes E_i [1]) && \longrightarrow \\
& \longrightarrow \Hom^{p + 1}_\D (\psi^\sharp_{E_i} (E_j), \psi^\sharp_{E_i} (E_{j'})) \longrightarrow \cdots. &&
\end{alignat*}
\normalsize
The identification $\psi^\sharp_{E_i} (E_j) \cong \R_{E_i} E_j$ implies that $\Hom^\bullet_\D (\psi^\sharp_{E_i} (E_j), E_i) \cong 0$.
Hence, we have
\[
\Hom^\bullet_\D (\psi^\sharp_{E_i} (E_j), \Hom^1_\D(E_{j'}, E_i)^* \otimes E_i [1]) \cong 
\Hom^1_\D(E_{j'}, E_i)^* \otimes \Hom^\bullet_\D (\psi^\sharp_{E_i} (E_j), E_i [1]) \cong 0.
\]
Therefore, there exists an isomorphism $\Hom^\bullet_\D (\psi^\sharp_{E_i} (E_j), \psi^\sharp_{E_i} (E_{j'})) \cong \Hom^\bullet_\D (\psi^\sharp_{E_i} (E_j), E_{j'})$.

Next, we see $\Hom^\bullet_\D (\psi^\sharp_{E_i} (E_j), E_{j'}) \cong \Hom^\bullet_\D (E_j, E_{j'})$.
By $\Hom^\bullet_\D (E_i, E_{j'}) \cong 0$, the long exact sequence 
\small
\begin{alignat*}{2}
& && \mathllap{\cdots} \longrightarrow \Hom^{p - 1}_\D (\psi^\sharp_{E_i} (E_j), E_{j'}) \longrightarrow \\
& \longrightarrow \Hom^p_\D (\Hom^1_\D(E_j, E_i)^* \otimes E_i [1], E_{j'}) \longrightarrow \Hom^p_\D (E_j, E_{j'}) && \longrightarrow \Hom^p_\D (\psi^\sharp_{E_i} (E_j), E_{j'}) \longrightarrow \\
& \longrightarrow \Hom^{p + 1}_\D (\Hom^1_\D(E_j, E_i)^* \otimes E_i [1], E_{j'}) \longrightarrow \cdots. &&
\end{alignat*}
\normalsize
yields an isomorphism $\Hom^\bullet_\D (\psi^\sharp_{E_i} (E_j), E_{j'}) \cong \Hom^\bullet_\D (E_j, E_{j'})$.
\qed
\end{pf}
\begin{lem}\label{lem : Hom 3 sharp}
Let $j = i^\sharp, \dots, i - 1$ and $l = i + 1, \dots, \mu$.
\begin{subequations}
\begin{enumerate}
\item If $\Hom^1_\D (E_j, E_i) \cong 0$ or $\Hom^\bullet_\D (E_i, E_l) \cong 0$, then we have 
\begin{equation}\label{eq : Hom 3 sharp 1}
\Hom^\bullet_\D (\psi^\sharp_{E_i} (E_j), E_l) \cong \Hom^\bullet_\D (E_j, E_l).
\end{equation}
\item If $\Hom^1_\D (E_j, E_i) \not\cong 0$, $\Hom^\bullet_\D (E_i, E_l) \not\cong 0$ and $\Hom^\bullet_\D (E_j, E_l) \cong 0$, then we have 
\begin{equation}\label{eq : Hom 3 sharp 2}
\Hom^p_\D (\psi^\sharp_{E_i} (E_j), E_l) \cong 
\begin{cases}
\Hom^1_\D (E_j, E_i) \otimes \Hom^p_\D (E_i, E_l), & p = p_{i, l}, \\
0, & p \ne p_{i, l}.
\end{cases}
\end{equation}
\item If $\Hom^1_\D (E_j, E_i) \not\cong 0$, $\Hom^\bullet_\D (E_i, E_l) \not\cong 0$ and $\Hom^\bullet_\D (E_j, E_l) \not\cong 0$, then we have 
\begin{equation}\label{eq : Hom 3 sharp 3}
\Hom^\bullet_\D (\psi^\sharp_{E_i} (E_j), E_l) \cong 0.
\end{equation}
\end{enumerate}
\end{subequations}
\end{lem}
\begin{pf}
If $\Hom^1_\D (E_j, E_i) \cong 0$, the exact triangle $\psi^\sharp_{E_i} (E_j) \to E_j \to \Hom^1_\D(E_j, E_i)^* \otimes E_i [1]$ implies $\psi^\sharp_{E_i} (E_j) \cong E_j$.
Hence, we have \eqref{eq : Hom 3 sharp 1}.

We assume $\Hom^1_\D (E_j, E_i) \not\cong 0$.
By the exact triangle $\psi^\sharp_{E_i} (E_j) \to E_j \to \Hom^1_\D(E_j, E_i)^* \otimes E_i [1]$, we have the long exact sequence
\small
\begin{align*}
\cdots & \overset{\epsilon_{p - 1}}{\longrightarrow} \Hom^{p - 1}_\D (E_j, E_l) \longrightarrow \Hom^{p - 1}_\D (\psi^\sharp_{E_i} (E_j), E_l) \longrightarrow \\
\longrightarrow \Hom^p_\D (\Hom^1_\D(E_j, E_i)^* \otimes E_i [1], E_l) & \overset{\epsilon_p}{\longrightarrow} \Hom^p_\D (E_j, E_l) \longrightarrow \Hom^p_\D (\psi^\sharp_{E_i} (E_j), E_l) \longrightarrow \\
\longrightarrow \Hom^{p + 1}_\D (\Hom^1_\D(E_j, E_i)^* \otimes E_i [1], E_l) & \overset{\epsilon_{p + 1}}{\longrightarrow} \Hom^{p + 1}_\D (E_j, E_l) \longrightarrow \cdots, 
\end{align*}
\normalsize
where $\epsilon_p \colon \Hom^p_\D (\Hom^1_\D(E_j, E_i)^* \otimes E_i [1], E_l) \longrightarrow \Hom^p_\D (E_j, E_l)$ denotes the $\CC$-linear map induced by the morphism $E_j \to \Hom^1_\D(E_j, E_i)^* \otimes E_i [1]$.
Note that we have 
\[
\Hom^{p}_\D (\Hom^1_\D(E_j, E_i)^* \otimes E_i [1], E_l) \cong \Hom^1_\D (E_j, E_i) \otimes \Hom^{p - 1}_\D (E_i, E_l).
\]
Hence, if we have $\Hom^\bullet_\D (E_i, E_l) \cong 0$, the long exact sequence implies the isomorphism \eqref{eq : Hom 3 sharp 1}.

We assume $\Hom^\bullet_\D (E_i, E_l) \not\cong 0$.
If $\Hom^\bullet_\D (E_j, E_l) \cong 0$, then the long exact sequence implies \eqref{eq : Hom 3 sharp 2}.

Next, we consider the last case $\Hom^\bullet_\D (E_j, E_l) \not\cong 0$.
Since the condition (E2) implies $p_{j, l} = p_{i, l} + 1$, the $\CC$-linear map $\epsilon_p$ is nontrivial for $p = p_{j, l}$.
By the condition (E1), the $\CC$-linear map $\epsilon_p$ can be identified with the map $\epsilon_p \colon \CC \longrightarrow \CC$.
We assume that $\epsilon_p$ is the zero map.
Then, by the long exact sequence, we have $\Hom^p_\D (\psi^\sharp_{E_i} (E_j), E_l) \cong \CC$ and $\Hom^{p - 1}_\D (\psi^\sharp_{E_i} (E_j), E_l) \cong \CC$.
However, this fact contradicts Proposition \ref{prop : strongly monochromaticity}.
Therefore, the $\CC$-linear map $\epsilon_p \colon \CC \longrightarrow \CC$ is an isomorphism and we obtain \eqref{eq : Hom 3 sharp 3}.
\qed
\end{pf}
Since the full exceptional collection $\E$ is Ext, we have 
\[
\Hom^p_\D (E_k, E_i [1]) \cong 
\begin{cases}
\Hom^{p}_\D (E_k, E_i), & p = p_{k, i} - 1, \\
0, & p \ne  p_{k, i} - 1,
\end{cases}
\]
for $k = 1, \dots, i^\sharp - 1$ and 
\[
\Hom^p_\D (E_i [1], E_l) \cong 
\begin{cases}
\Hom^{p}_\D (E_i, E_l), & p = p_{i, l} + 1, \\
0, & p \ne  p_{i, l} + 1,
\end{cases}
\]
for $l = i + 1, \dots, \mu$.
Note $p_{k, i} \ge 2$ by definition.
Therefore, the full exceptional collection $\E^\sharp_i$ is monochromatic and Ext.
Next, we see that the full exceptional collection $\E^\sharp_i$ satisfies the conditions (E1) and (E2).
Lemma \ref{lem : Hom 1 sharp}, \ref{lem : Hom 2 sharp}, \ref{lem : Hom 3 sharp} and \ref{lem : Hom 4 sharp} imply that $\E^\sharp_i$ satisfies the condition (E1).
\begin{lem}
The monochromatic full Ext-exceptional collection $\E^\sharp_i$ satisfies the condition (E2).
\end{lem}
\begin{pf}
For simplicity, we put $\E^\sharp_i = (E^\sharp_1, \dots, E^\sharp_\mu)$ and $p^\sharp_{a, b} = \deg \, \Hom^\bullet_\D (E^\sharp_a, E^\sharp_b)$.
By Lemma \ref{lem : Hom 1 sharp}, \ref{lem : Hom 2 sharp}, \ref{lem : Hom 3 sharp} and \ref{lem : Hom 4 sharp}, we have $\Hom^\bullet_\D (E^\sharp_a, E^\sharp_b) \not\cong 0$, $\Hom^\bullet_\D (E^\sharp_b, E^\sharp_c) \not\cong 0$ and $\Hom^\bullet_\D (E^\sharp_a, E^\sharp_c) \not\cong 0$ for $a < b < c$ if and only if the tuple $(E^\sharp_a, E^\sharp_b, E^\sharp_c)$ is one of the following:
\begin{enumerate}
\item $(E^\sharp_a, E^\sharp_b, E^\sharp_c) = (E_a, E_i [1], \psi^\sharp_{E_i} (E_{c - 1}))$ satisfying 
\[
\Hom^\bullet_\D (E_a, E_i) \not\cong 0, \quad \Hom^\bullet_\D (E_{c - 1}, E_i) \not\cong 0, \quad \Hom^\bullet_\D (E_a, E_{c - 1}) \cong 0.
\]
\item $(E^\sharp_a, E^\sharp_b, E^\sharp_c) = (E_i [1], \psi^\sharp_{E_i} (E_{b - 1}), E_c)$ satisfying 
\[
\Hom^\bullet_\D (E_i, E_{b - 1}) \not\cong 0, \quad \Hom^\bullet_\D (E_i, E_c) \not\cong 0, \quad \Hom^\bullet_\D (E_{b - 1}, E_c) \cong 0.
\]
\item $(E^\sharp_a, E^\sharp_b, E^\sharp_c) = (E_a, E_i [1], E_c)$ satisfying 
\[
\Hom^\bullet_\D (E_a, E_i) \not\cong 0, \quad \Hom^\bullet_\D (E_i, E_c) \not\cong 0, \quad \Hom^\bullet_\D (E_a, E_c) \not\cong 0.
\]
\item $(E^\sharp_a, E^\sharp_b, E^\sharp_c) = (E_a, E_b, E_c)$ satisfying 
\[
\Hom^\bullet_\D (E_a, E_b) \not\cong 0, \quad \Hom^\bullet_\D (E_b, E_c) \not\cong 0, \quad \Hom^\bullet_\D (E_a, E_c) \not\cong 0.
\]
\end{enumerate}
In the first case (i), by Lemma \ref{lem : Hom 2 sharp} we have 
\begin{align*}
p^\sharp_{a, b} & = \deg \, \Hom^\bullet_\D (E_a, E_i [1]) = p_{a, i} - 1, \\
p^\sharp_{b, c} & = \deg \, \Hom^\bullet_\D (E_i [1], \psi^\sharp_{E_i} (E_{c - 1})) = 1, \\
p^\sharp_{a, c} & = \deg \, \Hom^\bullet_\D (E_a, \psi^\sharp_{E_i} (E_{c - 1})) = p_{a, i}.
\end{align*}
In the second case (ii), by Lemma \ref{lem : Hom 3 sharp} we have 
\begin{align*}
p^\sharp_{a, b} & = \deg \, \Hom^\bullet_\D (E_i [1], \psi^\sharp_{E_i} (E_{b - 1})) = 1, \\
p^\sharp_{b, c} & = \deg \, \Hom^\bullet_\D (\psi^\sharp_{E_i} (E_{b - 1}), E_c) = p_{i, c}, \\
p^\sharp_{a, c} & = \deg \, \Hom^\bullet_\D (E_i [1], E_c) = p_{i, c} + 1.
\end{align*}
In the third case (iii), we have 
\begin{align*}
p^\sharp_{a, b} & = \deg \, \Hom^\bullet_\D (E_a, E_i [1]) = p_{a, i} - 1, \\
p^\sharp_{b, c} & = \deg \, \Hom^\bullet_\D (E_i [1], E_c) = p_{i, c} + 1, \\
p^\sharp_{a, c} & = \deg \, \Hom^\bullet_\D (E_a, E_c) = p_{a, c}.
\end{align*}
In the last case (iv), we have 
\begin{align*}
p^\sharp_{a, b} & = \deg \, \Hom^\bullet_\D (E_a, E_b) = p_{a, b}, \\
p^\sharp_{b, c} & = \deg \, \Hom^\bullet_\D (E_b, E_c) = p_{b, c}, \\
p^\sharp_{a, c} & = \deg \, \Hom^\bullet_\D (E_a, E_c) = p_{a, c}.
\end{align*}
Therefore, we obtain $p^\sharp_{a, b} + p^\sharp_{b, c} = p^\sharp_{a, c}$ by direct calculations.
\qed
\end{pf}
The last statements $\A^\sharp_{E_i} \cong \ex{\E^\sharp_i}$ and \eqref{eq : Sim E sharp} follows from Proposition \ref{prop : Exchange graph}.
We complete the proof of Proposition \ref{prop : tilting and full exceptional collection}.
\qed
\end{pf}
By Proposition \ref{prop : std heart} and \ref{prop : tilting and full exceptional collection}, we obtain the statement of Theorem \ref{thm : main thm}.
\qed
%%%%%%%%%%%%%%%%%%%%%%%%%%%%%%%%%%%%%%%%%%%%%%%%%%%%%%%%%%%%%%%%


\begin{thebibliography}{99}
{\small 

\bibitem[BM]{BM}
A.~Bayer, E.~Macr\`{i}, 
\textit{The space of stability conditions on the local projective plane},
Duke Math. J. {\bf 160} (2011), no. 2, 263–322.

\bibitem[BP]{BP} 
A.~Bondal and A.~Polishchuk, 
\textit{Homological Properties of Associative Algebras: The Method of Helices}, 
Izv. RAN. Ser. Mat., 1993, Volume 57, Issue 2, 3-50  (Mi izv877). 

\bibitem[B1]{B}
T.~Bridgeland,
\textit{Stability conditions on triangulated categories}, 
Ann. of Math. (2), {\bf 166} (2) : 317-345, 2007.

\bibitem[B2]{B2}
T.~Bridgeland, 
\textit{Spaces of stability conditions}, 
Algebraic geometry-Seattle 2005. Part 1, 1–21, Proc. Sympos. Pure Math., {\bf 80}, Part 1, Amer. Math. Soc., Providence, RI, 2009.

\bibitem[BQS]{BQS}
T.~Bridgeland, Y.~Qiu and T.~Sutherland,
\textit{Stability conditions and $A_2$-quiver}, 
Advances in Mathematics, Volume {\bf 365}, 13 May 2020, 107049. 

\bibitem[BBD]{BBD}
A.~Beilinson, J.~Bernstein and P.~Deligne,
\textit{Faisceaux Pervers, Ast\'{e}risque} {\bf 100},
Soc. Math de France, Paris (1983).

\bibitem[DK1]{DK}
G.~Dimitrov and L.~Katzarkov, 
\textit{Non-semistable exceptional objects in hereditary categories},
Int. Math. Res. Not. IMRN 2016, no. {\bf 20}, 6293–6377.

\bibitem[DK2]{DK2}
G.~Dimitrov and L.~Katzarkov, 
\textit{Non-semistable exceptional objects in hereditary categories: some remarks and conjectures},
Stacks and categories in geometry, topology, and algebra, 263–287, Contemp. Math., {\bf 643}, Amer. Math. Soc., Providence, RI, 2015.

\bibitem[DK3]{DK3}
G.~Dimitrov and L.~Katzarkov, 
\textit{Bridgeland stability conditions on the acyclic triangular quiver}, 
Adv. Math. {\bf 288} (2016), 825–886.

\bibitem[DK4]{DK4}
G.~Dimitrov and L.~Katzarkov, 
\textit{Bridgeland stability conditions on wild Kronecker quivers}, 
Adv. Math. {\bf 352} (2019), 27–55.

\bibitem[HKK]{HKK}
F.~Haiden, L.~Katzarkov, M.~Kontsevich,
\textit{Flat surfaces and stability structures},
Publ. Math. Inst. Hautes Études Sci. {\bf 126} (2017), 247–318.

\bibitem[KV]{KV}
B.~Keller and D.~Vossieck, 
\textit{Aisles in derived categories}, 
Bull. Soc. Math. Belg. {\bf 40} (1988), 239-253.

\bibitem[KQ]{KQ}
A.~King and Y.~Qiu,
\textit{Exchange graphs and Ext quivers},
Adv. Math. {\bf 285} (2015), 1106–1154.

\bibitem[L]{L}
C.~Li,
\textit{The space of stability conditions on the projective plane}, 
Selecta Math. (N.S.) {\bf 23} (2017), no. 4, 2927–2945.

\bibitem[M1]{M}
E.~Macr\`{i}, 
\textit{Stability conditions on curves},
Math. Res. Lett. {\bf 14} (2007), no. 4, 657–672.

\bibitem[M2]{M2}
E.~Macr\`{i}, 
\textit{Some examples of spaces of stability conditions on derived categories},
arXiv:math/0411613.

\bibitem[O]{O}
S.~Okada, 
\textit{Stability manifold of $\PP^1$}, 
J. Algebraic Geom. {\bf 15} (2006), no. 3, 487–505.

\bibitem[Q]{Q}
Y.~Qiu, 
\textit{Stability conditions and quantum dilogarithm identities for Dynkin quivers},
Adv. Math. {\bf 269} (2015), 220–264.

\bibitem[T]{T}
A.~Takahashi,
\textit{Matrix Factorizations and Representations of Quivers I},
arXiv:math/0506347.

\bibitem[W]{W}
J.~Woolf, 
\textit{Stability conditions, torsion theories and tilting},
J. Lond. Math. Soc. (2) {\bf 82} (2010), no. 3, 663–682.

}
\end{thebibliography}
\end{document}